\documentclass[12pt]{amsart}
\usepackage{amsmath,amssymb,amsthm,amscd}
\def\qed{\vbox{\hrule
\hbox{\vrule\hbox to 5pt{\vbox to 8pt{\vfil}\hfil}\vrule}\hrule}}

\newcommand{\beg}{\begin{eqnarray*}}
\newcommand{\begn}{\begin{eqnarray}}
\newcommand{\en}{\end{eqnarray*}}
\newcommand{\enn}{\end{eqnarray}}

\newtheorem{prop}{Proposition}[section]
\newtheorem{theo}[prop]{Theorem}
\newtheorem{lemm}[prop]{Lemma}

\newtheorem{rema}[prop]{Remark}

\newtheorem{defi}[prop]{Definition}

\begin{document}

\title{ Asymptotically Hyperbolic Metrics
on Unit Ball Admitting  Multiple Horizons }
\author{ ZhenYang Li}
\address{Key Laboratory of Pure and Applied mathematics, School of Mathematics Science, Peking University,
Beijing, 100871, P.R. China.} \email{lzymath@163.com}

 \author{ YuGuang Shi}
\address{Key Laboratory of Pure and Applied mathematics, School of Mathematics Science, Peking University,
Beijing, 100871, P.R. China.} \email{ygshi@math.pku.edu.cn}

 \author{ Peng Wu}
\address{Key Laboratory of Pure and Applied mathematics, School of Mathematics Science, Peking University,
Beijing, 100871, P.R. China.} \email{wupenguin@gmail.com}

\renewcommand{\subjclassname}{%
  \textup{2000} Mathematics Subject Classification}

\thanks{The research of  the second author is partially Supported by NKBRPC (2006CB805905)  and
Fok YingTong Education Foundation}

\begin{abstract}
In this paper, we construct an asymptotically hyperbolic metric with
scalar curvature -6 on unit ball $\mathbf{D}^3$, which contains
multiple horizons.
\end{abstract}

\maketitle \markboth{Zhenyang Li,  Yuguang Shi,  Peng Wu}
{Asymptotically Hyperbolic Metric on Unit Ball Admitting Multiple
Horizons }

\section{Introduction}
\setcounter{equation}{0} \hspace{0.4cm}

In general relativity,  the initial data set of Cauchy problem for
Einstein equations which is denoted by $(M,g_{ij},p_{ij})$ is of
great importance,  here $(M, g_{ij})$ is a complete Riemannian
3-manifold, $p_{ij}$ is a symmetric $2$-tensors on $M$ satisfying
constrain equations (see \cite{SY}). Among all of the initial data
sets, those with asymptotically flat(AF)(see \cite{SY}) and
asymptotically hyperbolic (AH)(see definition 2.1) metrics are of
most interests so far.

On the other hand,  horizon which is defined by a surface $\Sigma
\subset M$ satisfying $H_{\Sigma}=tr_{\Sigma}(p)$(see \cite{BC})
is very interesting geometric object. When $p=0$(i.e. time
symmetric case), the horizon is nothing but a minimal surface. The
Schwarzschild and anti-de Sitter-Schwarzschild space are the
simplest examples for AF and AH manifold with horizon
respectively. But, they both have non-trivial topology. However,
for physical and mathematical reason, people intend to construct
topologically trivial manifolds with horizons. In \cite{BM}, R.
Beig and N. \'O Murchadha show that there exists AF metric, which
contains horizon, with scalar flat on $\mathbb{R}^{3}$.  And Miao
in \cite{M} also construct the same kind AF manifolds by making
use of Schwarzschild metric and the conformal deformation. These
results offer examples of globally regular and AF initial data for
the Einstein vacuum equations with minimal surfaces.  Combining
the method of Miao with that of Chru\'sciel and Delay(\cite{CD}),
the author of \cite{C} gives an example of scalar flat AF metric
on $\mathbb{R}^{3}$ admitting multiple horizons.

In recent years, AH manifolds have drawn  more and more attention
of both mathematicians  and physicians. They arise when
considering solutions to the Einstein field equations with a
negative cosmological constant, or when considering "hyperboloidal
hypersurfaces" in space-times which are asymptotically flat in
isotropic directions. The horizons in AH manifolds are more than
minimal surfaces, since AH manifolds can be realized as an
asymptotically null   spacelike hypersurface in asymptotically
flat space-time. Therefore, in AH manifold, horizons refer not
only to boundaries of domains which are minimal surfaces(in the
case of considering  negative cosmological constant) but also to
boundaries satisfying $H=\pm2$ (in the case of asymptotically null
spacelike hypersurface in AF spacetimes). Recently, in \cite{ST},
the authors provide an example of AH manifold with constant scalar
curvature $-6$ and horizons(see also Theorem 2.2). Their main idea
is to glue the anti-de Sitter-Schwarzchild space with a ball and
deform it conformally for several times, then they get the desired
manifold. Furthermore, they can prove that the mass of their
example can be arbitrary large or small.  So, it is nature to
consider of constructing such an AH manifold admitting multiple
horizons. More precisely, in this paper we show that there exists
an AH metric on unit ball $\mathbf{D}^{3}$ with constant scalar
curvature $R=-6$ and multiple horizons(see Theorem 2.4). First, we
will construct a metric on unit 3-ball with multiple horizons
using cut-and-glue method. Secondly,  we will conformally deform
the metric to an AH metric with constant scalar curvature $R=-6$
by solving a nonlinear PDE, then, as in \cite{ST}, the existence
of horizons follows from implicit function theorem. We'd like to
remark that it seems that our method should work for the
construction of AF manifold with multiple horizons as what has
been done in \cite{C}.

The outline of this paper is as follows, In Sect.2 we  cut the
parts containing  horizons from some examples  of AH manifolds
given by \cite{ST} and glue them together smoothly. Consequently,
we obtain a new AH metric with multiple horizons. In Sect.3, we
perturb the  new metric conformally to an AH metric with constant
scalar curvature $-6$. Then by keeping the location of horizons
far enough to each other, we show that the existence of the
multiple horizons is guaranteed by a Lemma in \cite{ST}.

\section{Construction of Asymptotically Hyperbolic Metric On Unit Ball with Multiple Horizons}
\setcounter{equation}{0} \hspace{0.4cm}

In this section, we will complete the first step of the proof of the
main result, namely, we will construct an asymptotically hyperbolic
metric on $\mathbf{D}^{3}$ admitting multiple horizons by gluing
arguments, but the scalar curvature may not equal to $-6$. First of
all, Let us recall some basic definitions and facts.
\begin{defi}
A complete noncompact Riemannian manifold $(X^{3},g)$ is said to
be asymptotically hyperbolic if there is a compact manifold
$(\overline{X},\overline{g})$ with boundary $\partial X$ and a
smooth function $t$ on $\overline{X}$ such that the following are
true: \\
(i)$X=\overline{X}\setminus\partial X.$\\
(ii)$t=0$ on $\partial X$, and $t>0$ on $X$.\\
(iii)$\overline{g}=t^{2}g$ extends to be $C^{3}$ up to the
boundary.\\
(iv)$|dt|_{\overline{g}}=1$ at $\partial X$.\\
(v)Each component $\Sigma$ of $\partial X$ is the standard two
sphere $(\mathbb{S}^{2}, g_{0})$ and there is a collar
neighborhood of $\Sigma$ where
$$g=\sinh^{-2}t(dt^{2}+g_{t})$$
with
$$g_{t}=g_{0}+\frac{t^{3}}{3}h+O(t^{4})$$
where $h$ is a $C^{2}$ symmetric two tensor on $\mathbb{S}^{2}$.
\end{defi}

It is proved in \cite{W} that for an AH manifold $(X^{3},g)$ with
scalar curvature $R_{g}\geq-6$, the mass of an end of $X$
corresponding to a boundary component $\Sigma$ of $\partial X$ is
well-defined and given by
$$M=\frac{1}{16\pi}[(\int_{\mathbb{S}^{2}}tr_{g_{0}}(h)dv_{g_{0}})^{2}-
|\int_{\mathbb{S}^{2}}tr_{g_{0}}(h)(x)xdv_{g_{0}}|^{2}]^{\frac{1}{2}}$$
where $x$ is the standard coordinates of a point on
$\mathbb{S}^{2}$ in $\mathbb{R}^{3}$.

We denote the standard hyperbolic space by $\mathbb{H}^{3}$, and
introduce the ball model for $\mathbb{H}^{3}$ which is denoted by
$(\mathbf{D^{3}},ds^{2}_{\mathbb{H}^{3}}) $, here, $\mathbf{D^{3}}$
is the unit ball in $\mathbb{R}^{3}$, and $ds^{2}_{\mathbb{H}^{3}}$
is the standard hyperbolic metric which is defined as following:

$$ds^{2}_{\mathbb{H}^{3}}=\frac{4}{(1-|x|^{2})^{2}}\displaystyle\sum_{i=1}^{3}(dx^{i})^{2},$$
where $\displaystyle\sum_{i=1}^{3}(dx^{i})^{2}$ is the Euclidean
metric.

In \cite{ST}, the  authors construct a family of asymptotically
hyperbolic metrics on
 $\mathbf{D}^{3}$ as follows:
\begin{theo}(\cite{ST})
Let $\mathbf{D}^{3}$ be the unit ball in $\mathbb{R}^{3}$. For any
$M>0$ and $\delta$, there is a smooth complete metric $g$ on
$\mathbf{D}^{3}$ with constant scalar curvature $-6$ such that the
following are true\\
(i)$(\mathbf{D}^{3}, g)$ is asymptotically hyperbolic with mass
$M_{g}$ satisfying $|M_{g}-M|<\delta$.\\
(ii)There exist surfaces $S_{1}, S_{2}$ and $S_{3}$ which are
topological spheres with constant mean curvature $-2, 0, 2$
respectively such that $S_{1}$ is in the interior of $S_{2}$ and
$S_{2}$ is in the interior of $S_{3}$.\\
(iii)Outside a compact set $U$, the metric g is conformal to the
standard hyperbolic metric of $\mathbf{D}^{3}$ , and $S_{1}, S_{2}$
and $S_{3}$ are contained in $\mathbf{D}^{3}\setminus U$.
\end{theo}

Let us fix some notations  for our paper. We will denote the
origin by $o$. Let $x\in\mathbf{D}^{3}$ and $B_{x}(\rho)$ be the
geodesic ball centered at $x$ with radius $\rho$ under the
standard hyperbolic metric on $\mathbf{D}^{3}$. The hyperbolic
distance starting from $x$ to $y$ will be denoted by $\rho_{x}(y)$
(simply by $\rho(y)$ if $x$ is the origin). Without loss of
generality, we may reformulate Theorem 2.2 in the following way:

\begin{theo}
Let $S_{1},S_{2},S_{3}$  be the surfaces and  $g$ be the AH metric
which are given by Theorem 2.2. Then there exists a geodesic ball
$B_{o}(\delta)$ under the hyperbolic metric, such that
$S_{1},S_{2},S_{3}$  are contained in $B_{o}(\delta)\setminus
B_{o}(\frac{\delta}{2})$ and $g$  is conformal to the standard
hyperbolic metric of $\mathbf{D}^{3}$ on $\mathbf{D}^{3}\setminus
B_{o}(\frac{\delta}{4})$.
\end{theo}

Using gluing method and conformal deformation again, we are able
to prove our main result:

\begin{theo}
Let $\mathbf{D}^{3}$ be the unit open  ball in $\mathbb{R}^{3}$.
For any $K>0$, there is an AH   metric $g$ on $\mathbf{D}^{3}$
with constant scalar curvature $-6$,  such that there is
$\{x_{k}\}_{k=1}^{K}\subset\mathbf{D}^{3} $ and surfaces
$S_{1}^{i}, S_{2}^{i}$ and $S_{3}^{i}$, $1\leq i\leq K$,  which
are topological spheres with constant mean curvature $-2, 0, 2$
respectively and contained in  $B_{x_{i}}(\delta)$
 which does not intersect each other;
moreover, $S_{1}^{i}$ is in the interior of $S_{2}^{i}$ and
$S_{2}^{i}$ is in the interior of $S_{3}^{i}$,  and outside a
compact set, the metric g is conformal to the standard hyperbolic
metric of $\mathbf{D}^{3}$.
\end{theo}

We only consider the case for $K=2$, since the other cases are
essentially the  same. By gluing arguments, we will show

\begin{prop} There is smooth AH  metric $\widetilde{g}$ on $\mathbf{D}^3$ with the
following properties:
\begin{enumerate}
  \item $B_{o}(\delta_{2})\setminus B_{o}(\frac{\delta_{2}}{2})$ and
$B_{p}(\delta_{1})\setminus B_{p}(\frac{\delta_{1}}{2})$ each
contains the surfaces  with mean curvature 2, 0, -2, here, $o$,
$p$ are two points in $\mathbf{D}^3$, with hyperbolic distance
being $2\tau\triangleq  \rho(p)>10 (\delta_1 +\delta_2)$, so that
$B_{o}(\delta_{2})$ and $B_{p}(\delta_{1})$ does not interest each
other. Without loss of generality, we may assume $\tau\geq 100$.

\item The scalar curvature $R_{\widetilde{g}}$ of $\widetilde{g}$ satisfying
$$
R_{\widetilde{g}}(x)=-6,
$$
for $x\in\mathbf{D}^{3}\setminus(B_{p}(\tau+2)\setminus
B_{p}(\tau+1))$ and
$$
|R_{\widetilde{g}}(x)+6|\leq Ce^{-3\tau}, $$ for $x\in
B_{p}(\tau+2)\setminus B_{p}(\tau+1)$, here, $C$ is a positive
constant which is independent of $\tau$.
\item $\widetilde{g}$ is  conformal to
the hyperbolic metric outside $ B_{o}(\frac{\delta_{2}}{4})\cup
B_{p}(\frac{\delta_{1}}{4}) $
\end{enumerate}
\end{prop}

\begin{proof}

Choose two AH metrics $g_{1}, g_{2}$ as given by Theorem 2.2(not
necessary having the same mass). Let $B_{o}(\delta_{1})$ and
$B_{o}(\delta_{2})$ be the sets described in Theorem 2.3 for
$(\mathbf{D}^{3},g_{1})$ and $ (\mathbf{D}^{3},g_{2})$ respectively,
such that for $i=1$, $2$,
\begin{eqnarray}
g_{i}(y)= \phi_{i}^{4}(y)ds^{2}_{\mathbb{H}^{3}},\nonumber
\end{eqnarray}
where $y\in\mathbf{D}^{3}\setminus B_{o}(\frac{\delta_{i}}{4})$.
And by Theorem 4.1 in \cite{ST}, $\phi_{i}$ satisfies:

\begin{eqnarray}
\parallel\phi_{i}(y)-1\parallel_{\mathbf{C}^{3}}\leq
Ce^{-3\rho(y)},
\end{eqnarray}
  for
$y\in\mathbf{D}^{3}\setminus B_{o}(\frac{\delta_{i}}{4})$  and $C$
is a positive constant that is independent of $y$.

Now, we will glue $g_{1}$ and $g_{2}$ together as following:

Let us introduce the upper halfspace model $\mathbb{R}^{3}_{+}$
for $\mathbb{H}^{3}$ and label the point $x\in\mathbb{H}^{3}$ by
$(\mathbf{x}, y)$ with $\mathbf{x}\in \mathbb{R}^{2} $  and
$y\in\mathbb{R}_{+}$. Under this coordinate system, the standard
metric for hyperbolic space can be expressed as:
$$ds^{2}_{\mathbb{H}^{3}}=\frac{(dx^{1})^{2}+(dx^{2})^{2}+dy^{2}}{y^{2}}$$
here, $\mathbf{x}=(x^{1}, x^{2})\in\mathbb{R}^{2}$.  Suppose
$o=(\mathbf{0},1)$ and $p=(\mathbf{x}_{p}, y_{p})$, then there is
a hyperbolic translation $F$ which maps $o$ to $p$ with
$\rho(p)\triangleq 2\tau >10(\delta_1+\delta_2)$
\begin{eqnarray}\begin{array}{ccc}
F: B_{o}(\tau+3)\longrightarrow B_{p}(\tau+3)\\
   \hspace{4cm}(\mathbf{x},y)\longmapsto F(\mathbf{x},y)= (\mathbf{x}_{p}+y_{p}\mathbf{x},
   y_{p}y).
\end{array}\nonumber\end{eqnarray}
Then  it is easy to see that $F$ induced a natural isometry
between the standard hyperbolic metric in  $B_p(\tau+3)$ and that
in $B_o(\tau+3)$:\\
\begin{eqnarray}&(B_p(\tau+3), (F^{-1})^*ds^{2}_{\mathbb{H}^{3}}|_{B_0(\tau+3)})
\cong (B_p(\tau+3), ds^{2}_{\mathbb{H}^{3}}|_{B_p(\tau+3)}),\\
&\nonumber\end{eqnarray} indeed,
$(F^{-1})^*(ds^{2}_{\mathbb{H}^{3}}|_{B_0(\tau+3)})=ds^{2}_{\mathbb{H}^{3}}|_{B_p(\tau+3)}$,
in the sense  that the metrics on the both sides have the same
components under the standard upper half space coordinates.

Therefore, we can pull $g_{1}$  on $B_o(\tau+3)$ to $B_p(\tau+3)$ by
the diffeomorphism $F$, which gives an isometry:
\begin{eqnarray}(B_p(\tau+3), (F^{-1})^*g_1) \cong (B_o(\tau+3), g_1).\nonumber\end{eqnarray}
Because of (2.2), we can identify $B_o(\tau+3)$ and $B_p(\tau+3)$
both equipped with hyperbolic metric via $F$, and
$$(F^{-1})^*g_1 = (\phi_1\circ F^{-1})^{4}(F^{-1})^*
(ds^{2}_{\mathbb{H}^{3}}|_{B_0(\tau+3)}) = (\phi_1\circ
F^{-1})^{4}ds^{2}_{\mathbb{H}^{3}}|_{B_p(\tau+3)},$$ for $x\in
B_p(\tau+3)\setminus B_p(\frac{\delta_{1}}{4})$. Also, the
inequality (2.1) for $\phi_{1}$ can be described as:
\begin{eqnarray}
\parallel\phi_{1}\circ F^{-1}(y)-1\parallel_{\mathbf{C}^{3}}\leq
Ce^{-3\rho_{p}(y)}, \nonumber\end{eqnarray} for $y\in
B_p(\tau+3)\setminus B_p(\frac{\delta_{1}}{4})$. For simplicity,
$\phi_{1}\circ F^{-1}(y)$ will still be denoted by $\phi_1 (y)$ in
the sequel.

Let $\eta$ be a smooth cut-off function such that $0\leq\eta\leq1$
and
\begin{eqnarray} \eta(x)=\left \{
         \begin{array}{lll}
            1&x\in B_{p}(\tau+1),\\\\
              0&x\in\mathbf{D}^{3}\setminus
              B_{p}(\tau+2).
\end{array}  \right.\nonumber
\end{eqnarray}
Hence $||\eta||_{C^{2}}$ is uniformly bounded. Next, we define a
new metric  $\widetilde{g}$ on $\mathbf{D}^{3}$, which is given by
\begin{equation}
\widetilde{g}(x)=\left\{
  \begin{array}{ll}
   (F^{-1})^*g_1, & \hbox{$x\in B_{p}(\tau+1)$;} \\
   (\eta\phi_1+(1-\eta)\phi_2)^4ds^{2}_{\mathbb{H}^{3}}|_{B_p(\tau+2)\setminus
B_p(\tau+1) }, & \hbox{$x\in B_p(\tau+2)\setminus
              B_{p}(\tau+1)$;} \\
    g_{2}, & \hbox{$x\in \mathbf{D}^{3}\setminus B_{p}(\tau+2)$.}
  \end{array}
\right.\nonumber
\end{equation}
By its definition, we see $\widetilde{g}$ satisfies (1) in
Proposition 2.5.

Again  by the definition of $\widetilde{g}$ and (2.1), we can
calculate that the scalar curvature $R_{\widetilde{g}}$ of
$\widetilde{g}$ satisfies that:
\begin{eqnarray}
R_{\widetilde{g}}(x)=-6\hspace{0.4cm}
\mbox{for}\hspace{0.4cm}x\in\mathbf{D}^{3}\setminus(B_{p}(\tau+2)\setminus
B_{p}(\tau+1))\nonumber
\end{eqnarray}
and for $x\in B_{p}(\tau+2)\setminus B_{p}(\tau+1)$,
\begin{eqnarray}
|R_{\widetilde{g}}(x)+6|\leq Ce^{-3\tau}.\nonumber
\end{eqnarray}
thus, we verified (2) in Proposition; and by Theorem 2.3, we see
that (3) in Proposition 2.5 is also true, therefore, we finish to
prove the Proposition.
\end{proof}
\begin{rema}One can see from the construction of $\widetilde{g}$
that it  depends on $\tau$. To emphasize this, we will denote
$\widetilde{g}$ by $\widetilde{g}_{\tau}$ in the next section.
\end{rema}

\section{Proof of the Main result by Conformal Deformation}
\setcounter{equation}{0} In this section, we will prove our main
result Theorem 2.4. Namely, we perturb $\widetilde{g}_{\tau}$
constructed in last section by conformal deformation, and show
that the resulting metric is an AH metric with scalar curvature
equal to $-6$ and containing multiple horizons. For this purpose,
we need

\begin{lemm} Let $\widetilde{g}_\tau$ be constructed in the Proposition 2.5,
then there is $u_\tau>0$ such that $g_{\tau}=u_\tau ^4
\widetilde{g}_\tau$ is an AH metric with scalar curvature $R=-6$,
and
$$\lim_{\tau\rightarrow \infty}\sup_{\mathbf{D}^3}|u_\tau-1| =0.
$$
Moreover,  outside a compact set, the metric $\widetilde{g}_\tau$
is conformal to the standard hyperbolic metric on
$\mathbf{D}^{3}$.
\end{lemm}

\begin{proof}

It is sufficient to solve the following equation:
\begin{eqnarray}
\left \{
         \begin{array}{lll}
            \triangle_{\widetilde{g}_\tau}u-\frac{3}{4}u^{5}-\frac{1}{8}R_{\widetilde{g}_\tau}u=0,\\\\
              \displaystyle\lim_{\rho(x)\rightarrow+\infty}u(x)=1.
\end{array}  \right.
\end{eqnarray}

To do this, we will use exhausting domain arguments,   let us
choose a sequence $\{\rho_{k}\}_{k=1}^{\infty}$ with
$\rho_{1}\geq3\tau$ such that $\rho_{k}<\rho_{k+1}$ for
$k\in\mathbf{N}$ and $\rho_{k}\rightarrow+\infty$ as
$k\rightarrow+\infty$. Consider the following Dirichlet problem:
\begin{eqnarray}
\left \{
         \begin{array}{lll}
            \triangle_{\widetilde{g}_\tau}u_{i}-\frac{3}{4}u_{i}^{5}-\frac{1}{8}R_{\widetilde{g}_\tau}u_{i}=0,
            \hspace{0.4cm}\mbox{on} \hspace{0.4cm}B_{o}(\rho_{i})\\\\
             u_{i}|_{\partial B_{o}(\rho_{i})}=\phi^{-1}_{2}.
\end{array}  \right.
\end{eqnarray}
By the standard variational method, we see (3.2) has a smooth
nonnegative solution. And by  maximal principle, $u_{i}$ must be
positive. By assumption for $\phi_{2}$, we know that
$1-Ce^{-3\rho_{i}}\leq u_{i}|_{\partial B_{o}(\rho_{i})}\leq1 +
Ce^{-3\rho_{i}}$. We claim that
\begin{equation}
\sup_{B_{o}(\rho_{i})}|u_i -1|\leq C e^{-3\tau},
\end{equation}
here, $C$ is independent of $i$ and $\tau$. Let us prove the lower
bounded estimate is true. Indeed, suppose $u$ attains its minimum at
$x_{0}$, if $x_{0}\in \partial B_{o}(\rho_{i})$, then the claim of
the lower bound follows, otherwise, at an interior point $x_0$, one
has

\begin{eqnarray}
0\leq u_{i}^{-1}\triangle_{\widetilde{g}_\tau}u_{i} & = &
\frac{3}{4}u_{i}^{4}+\frac{1}{8}R_{\widetilde{g}_{\tau}} \nonumber
\end{eqnarray}
By (2) of Proposition 2.5, we get the claim for the lower bound, by
the similar arguments, we will get the upper bounded estimate. Thus
claim is true. And then, by exhausting domain, we get the solution
$u_\tau$  of (3.1), and

$$\sup_{\mathbf{D}^3}|u_\tau -1|\leq C e^{-3\tau},
$$
here, $C$ is a constant which is  independent of $\tau$.

Next, we will construct barrier of the equation  at the infinity
of the manifold, and by using this we show that $u$ approaches to
$1$ at desired rate.

 Let $v_{i}=u_{i}\phi_{2}$.
Since $g_{2}$ is  conformal to the standard hyperbolic metric
outside $B_{o}(3\tau)$, we have
\begin{eqnarray}
\hspace{1cm}\left \{
         \begin{array}{lll}
            \mathbf{L}(v_i):=\triangle_{\mathbb{H}^{3}}v_{i}-\frac{3}{4}v_{i}(v_{i}^{4}-1)=0,
            \hspace{0.4cm}\mbox{for}\hspace{0.4cm} x\in B_{o}(\rho_{i})\setminus
B_{o}(\overline{\rho})\\\\
              v_{i}|_{\partial B_{o}(\rho_{i})}=1,
\end{array}  \right.
\nonumber\end{eqnarray}
 for some $\overline{\rho}\geq3\tau$ where $\triangle_{\mathbb{H}^{3}}$ is the Laplacian operator for
hyperbolic metric.  Set $f_{-}(x)=1-\lambda e^{-3\rho(x)}$, for
$\lambda\geq0$. Hence,
\begin{eqnarray}
\mathbf{L}(f_{-})&=&-9\lambda e^{-3\rho(x)}+6\lambda e^{-3\rho(x)}\coth\rho\nonumber\\
&&-\frac{3}{4}(1-\lambda  e^{-3\rho(x)})[
            (1-\lambda  e^{-3\rho(x)})^{4}-1]\nonumber\\
            &=&6\lambda(\coth\rho-1) e^{-3\rho(x)}+O( e^{-6\rho(x)}).\nonumber
\end{eqnarray}
 Since
$(\coth\rho-1) e^{-3\rho(x)}>2 e^{-5\rho(x)}$, hence,   for
sufficiently large $\overline{\rho}$ we have $L(f_{-})>0$ in
$B_o(\rho_i)\setminus B_o(\bar \rho)$. So let us choose
$\lambda=e^{3\overline{\rho}}$,  then $\mathbf{L}(f_{-})>0$
whenever $\rho\geq\overline{\rho}$, $v_i =1\geq f_{-}$ at
$\partial B_{o}(\rho_{i})$, and  $v_i \geq
 f_{-}=0$ at $\partial B_{o}(\overline{\rho})$. Due to (2.1) and
 (3.3), by choosing $\tau$ and $\bar\rho$ sufficiently large, we
 may assume  $v_i \geq 2^{-\frac12}>5^{-\frac14}$ on $B_{o}(\rho_{i})\setminus
B_{o}(\overline{\rho})$. Now, we claim that $v_i \geq f_{-}$ on
$B_{o}(\rho_{i})\setminus B_{o}(\overline{\rho})$. For any point of
the boundary, the claim is obviously true. Suppose the claim fails,
then for any $p\in  B_{o}(\rho_{i})\setminus B_{o}(\overline{\rho})$
 with $(v_i -f_{-})(p)= \inf_{B_{o}(\rho_{i})\setminus
B_{o}(\overline{\rho})}(v_i -f_{-})<0$ is an interior point, and at
which we have

$$5^{-\frac14}<2^{-\frac12}\leq v_i (p)< f_{-}(p),$$
thus, at $p$, we have
$$0\leq \triangle_{\mathbb{H}^{3}}(v_i -f_{-}) \leq \frac34 v_i (v_i
^4 -1)- \frac34f_{-}(f_{-}^4 -1) <0,
$$
which is a contradiction. Therefore, we obtain

 \begin{eqnarray}v_{i}-1\geq -\lambda e^{-3\rho},
  \hspace{0.4cm}\mbox{in}\hspace{0.4cm} B_{o}(\rho_{i})\setminus B_{o}(\overline{\rho}).\nonumber
  \end{eqnarray}
On the other hand, by noting supper solution $f_{+}(x)=1+\lambda
e^{-3\rho(x)}$ and using the similar arguments we can  prove that
$v_{i}-1\leq \lambda e^{-3\rho}$ in $B_{o}(\rho_{i})\setminus
B_{o}(\overline{\rho})$. Hence, we have
$$|v_{i}-1|\leq \lambda e^{-3\rho(x)},\hspace{0.4cm}\mbox{in}\hspace{0.4cm} B_{o}(\rho_{i})\setminus B_{o}(\overline{\rho}).$$
Consequently, $u_\tau$ satisfies that
\begin{eqnarray}|u_\tau -1|\leq \lambda e^{-3\rho}\hspace{0.4cm}\mbox{in}\hspace{0.4cm}
\mathbf{D}^{3}\setminus B_{o}(\overline{\rho}).
\end{eqnarray}
and also by (3.3)
\begin{eqnarray}\sup_{\mathbf{D}^3}|u_\tau -1|\leq Ce^{-3\tau}\end{eqnarray}
as mentioned above here $C$ is a constant that is independent of
$\tau$.

Now applying Lemma 4.2 in \cite{ST},  we conclude from (3.4) that
the manifold $(\mathbf{D^{3}},g_{\tau})$ with
$g_{\tau}=u_{\tau}^{4}\widetilde{g}_\tau$ is an AH manifold. And
it follows from Proposition 2.5 that $g_{\tau}$ is conformal to
the hyperbolic metric outside some compact subset of
$\mathbf{D}^{3}$. Thus, we finish to prove the Lemma.
\end{proof}

Next, we have the following
\begin{lemm} Let $g_{\tau}$ be as in Lemma 3.1,which depends on $\tau$.
Then for sufficiently large $\tau$,
 $(\mathbf{D}^3, g_{\tau})$ contains  surfaces $S_{1}^{i},
S_{2}^{i}$ and $S_{3}^{i}$, $1\leq i\leq 2$,  which are
topological spheres with constant mean curvature $-2, 0, 2$
 and contained in $B_{o}(\delta)$  and $B_{p}(\delta)$ respectively which does not
intersect each other; moreover, $S_{1}^{i}$ is in the interior of
$S_{2}^{i}$ and $S_{2}^{i}$ is in the interior of $S_{3}^{i}$.
\end{lemm}
\begin{proof} Let us show the horizons  are in $(B_o (2\delta), g_{\tau})$
when $\tau$ is large enough. In fact, by the Lemma 4.4 in \cite{ST},
it is sufficient to show that
\begin{equation}
\|u_{\tau}-1\|_{C^{2, \alpha}(B_o (2\delta)\setminus B_o
(\frac{\delta}{4}))}\leq \epsilon,
\end{equation}
here $\epsilon$ and $u_{\tau}$ are given in Lemma 4.4 in \cite{ST}
and Lemma 3.1 respectively. Note in $(B_o (2\delta),
\widetilde{g}_\tau)$, sectional curvature is bounded, injective
radius has a uniform positive lower bound, then $B_o (\frac32
\delta)$ can be covered by finite number of harmonic coordinates
which has uniform size, (for existence of harmonic coordinates and
estimates of its size please see \cite{JK}). The number  and the
size of these harmonic coordinates are independent of $\tau$. For
any $x\in B_o (\frac32 \delta)$, without loss of its generality,
we may assume $B_x (1)$ has already covered by such harmonic
coordinates, then the components of the metric $g_{\tau}$ under
the harmonic coordinates satisfies

$$\parallel (\widetilde{g}_\tau)_{ij}\parallel_{C^{1, \alpha}}\leq C,$$
here $C$ is independent of $\tau$ and $x$. Now in $B_x (1)$, we have
equation

$$\triangle_{\widetilde{g}_\tau}u_{\tau}=\frac{3}{4}u_{\tau}(u_{\tau}^{4}-1),$$
by (3.1). Then combining (3.5) with  the standard estimate of PDE,
we get
$$\|u_{\tau}-1\|_{C^{2, \alpha}(B_x (1))}\leq C e^{-3\tau},
$$
where $C$ is independent of $\tau$. Since $x$ is arbitrary in $B_o
(2\delta)$, and $\tau$ can be arbitrary large, we get (3.6), which
implies there are horizons in $B_o (2\delta)$, by the same
arguments, we can show there is also horizons in $B_p (2\delta)$.
Thus we finish to prove the Lemma.
\end{proof}

Combine the above Lemmas, we get the proof of Theorem 2.4, which
finish to prove our main result.

\bibliographystyle{amsplain}

\end{document}